\documentclass[11pt]{article}
\usepackage{amsmath,amsthm,amsfonts}
\newcommand\pd[2]{\frac{\partial#1}{\partial#2}}
\renewcommand{\=}{\doteq}

\newtheorem{thm}{Theorem}[section]
 \newtheorem{prop}[thm]{Proposition}
 \newtheorem{lemma}[thm]{Lemma}
\theoremstyle{definition}

\theoremstyle{definition}
 \newtheorem{rem}[thm]{Remark}
\numberwithin{equation}{section}
\numberwithin{equation}{section}
\begin{document}
\title{\bf  Moufang symmetry IV.\\
Reductivity and hidden associativity
}
\author{Eugen Paal}
\date{}
\maketitle
\thispagestyle{empty}
\begin{abstract}
It is shown how integrability of the generalized Lie equations of a local analytic Moufang loop is related to the reductivity conditions and Sagle-Yamaguti identity. 
\par\smallskip
{\bf 2000 MSC:} 20N05, 17D10
\end{abstract}

\section{Introduction}

In this paper we proceed explaing the Moufang symmetry. It is shown how integrability of the generalized Lie equations of a local analytic Moufang loop is related to the reductivity conditions and Sagle-Yamaguti identity. The paper can be seen as a continuation of \cite{Paal1,Paal2,Paal3}.

\section{Generalized Lie equations}

In \cite{Paal1} the \emph{generalized Lie equations} (GLE) of a local analytic Moufang loop $G$ were found. These read
\begin{subequations}
\label{gle}
\begin{align}
w^{s}_{j}(g)\pd{(gh)^{i}}{g^{s}}+u^{s}_{j}(h)\pd{(gh)^{i}}{h^{s}}+u^{i}_{j}(gh)&=0\\
v^{s}_{j}(g)\pd{(gh)^{i}}{g^{s}}+w^{s}_{j}(h)\pd{(gh)^{i}}{h^{s}}+v^{i}_{j}(gh)&=0\\
u^{s}_{j}(g)\pd{(gh)^{i}}{g^{s}}+v^{s}_{j}(h)\pd{(gh)^{i}}{h^{s}}+w^{i}_{j}(gh)&=0
\end{align}
\end{subequations}
where $gh$ is the product of $g$ and $h$, and the auxiliary functions $u^s_j$, $v^s_j$ and $w^s_j$ are related with the  constraint
\begin{equation}
u^s_j(g)+v^s_j(g)+w^s_j(g)=0
\end{equation}
In this paper we inquire integrability of GLE (\ref{gle}a--c). Triality \cite{Paal2} considerations are wery helpful.

\section{Generalized Maurer-Cartan equations and Yamagutian}

Recall from \cite{Paal1} that for $x$ in  $T_e(G)$ the infinitesimal translations of $G$ are defined by
\begin{equation*}
L_x\=x^j u^s_j(g)\pd{}{g^s},\quad
R_x\=x^j v^s_j(g)\pd{}{g^s},\quad
M_x\=x^j w^s_j(g)\pd{}{g^s}\quad \in T_g(G)
\end{equation*}
with constriant
\begin{equation*}
L_x+R_x+M_x=0
\end{equation*}
Following triality \cite{Paal2} define the Yamagutian $Y(x;y)$ by
\begin{equation*}
6Y(x;y)=[L_x,L_y]+[R_x,R_y]+[M_x,M_y]
\end{equation*}
We know  from \cite{Paal2} the generalized Maurer-Cartan equations:
\begin{subequations}
\label{m-c}
\begin{align} 
[L_{x},L_{y}]&=L_{[x,y]}-2[L_{x},R_{y}]\\
[R_{x},R_{y}]&=R_{[y,x]}-2[R_{x},L_{y}]\\
[L_{x},R_{y}]&=[R_{x},L_{y}]
\end{align}
\end{subequations}
The latter can be rewritten \cite{Paal2} as follows:
\begin{subequations}
\label{lr-yam}
\begin{align}
[L_{x},L_{y}]=2Y(x;y)+\frac{1}{3}L_{[x,y]}+\frac{2}{3}R_{[x,y]}\\
[L_{x},R_{y}]=-Y(x;y)+\frac{1}{3}L_{[x,y]}-\frac{1}{3}R_{[x,y]}\\
[R_{x},R_{y}]=2Y(x;y)-\frac{2}{3}L_{[x,y]}-\frac{1}{3}R_{[x,y]}
\end{align}
\end{subequations}
Nonassociativity can be described \cite{Paal1} in $G$ by commutator 
\begin{equation*}
[L_x,R_y]=x^jy^k l^s_{jk}(g)\pd{}{g^s}
\end{equation*}
with the second-order associator coefficients 
\begin{equation*}
l^s_{jk}
\=
v^{s}_{j}(g)\pd{u^{i}_{k}(g)}{g^{s}}
               -u^{s}_{k}(g)\pd{v^{i}_{j}(g)}{g^{s}}
\end{equation*}

\section{Reductivity}

Define the  (secondary) auxiliary functions of $G$ by
\begin{align*}
u^s_{jk}(g)
&\=u^p_k(g)\pd{u^s_j(g)}{g^p}-u^p_j(g)\pd{u^s_k(g)}{g^p}\\
v^s_{jk}(g)
&\=v^p_k(g)\pd{v^s_j(g)}{g^p}-v^p_j(g)\pd{v^s_k(g)}{g^p}\\
w^s_{jk}(g)
&\=w^p_k(g)\pd{w^s_j(g)}{g^p}-w^p_j(g)\pd{w^s_k(g)}{g^p}
\end{align*}
The Yamaguti functions $Y^i_{jk}$ are defined by
\begin{equation*}
6Y^s_{jk}(g)\=u^s_{jk}(g)+v^s_{jk}(g)+w^s_{jk}(g)
\end{equation*}
In \cite{Paal3} we proved

\begin{thm}
The integrability conditons of the GLE (\ref{gle}a--c) read
\begin{equation}
\label{gle2yam}
Y^s_{jk}(g)\pd{(gh)^i}{g^s}+Y^s_{jk}(h)\pd{(gh)^i}{g^s}=Y^i_{jk}(gh)
\end{equation}
\end{thm}
Consider the first-order approximation of the integrability conditions (\ref{gle2yam}). We need 
\begin{lemma}
One has
\begin{equation}
\label{assoc2yam}
Y^i_{jk}=l^i_{jk}+\frac{1}{3}C^s_{jk}(u^i_s-v^i_s)
\end{equation}
\end{lemma}

\begin{proof}
Use formula (\ref{lr-yam}b).
\end{proof}

Introduce the Yamaguti constants $Y^i_{jkl}$ by
\begin{equation*}
Y^i_{jk}(g)=Y^i_{jkl}g^l+O(g^2)
\end{equation*} 
Then, by defining \cite{Paal1} the third-order associators $l^i_{jkl}$ by 
\begin{equation*}
l^i_{jk}(g)=l^i_{jkl}g^l+O(g^2)
\end{equation*}
it follows from Lemma \ref{assoc2yam} that
\begin{equation}
\label{yam3}
Y^i_{jkl}=l^i_{jkl}+\frac{1}{3}C^s_{jk}C^i_{sl}
\end{equation}
Now we can calculate:
\begin{gather*}
Y^i_{jk}(gh)
=Y^i_{jk}(h)+\pd{Y^i_{jk}(h)}{h^s}u^s_l(h)g^l+O(g^2)\\
Y^i_{jk}(g)\pd{(gh)^i}{g^s}+Y^s_{jk}(h)\pd{(gh)^i}{h^s}
=Y^i_{jkl}g^l u^i_s(h)+Y^s_{jk}(h)\left(\delta^i_s+\pd{u^i_{l}(h)}{h^s}g^l\right)+O(g^2)
\end{gather*}
By substituting the latter into (\ref{gle2yam}), comparing the coefficients at $g^l$ and replacing $h$ by $g$  we obtain the \emph{reductiovity conditions}
\begin{equation}
\label{pre-red}
u^s_l(g)\pd{Y^i_{jk}(g)}{h^s}-Y^i_{jk}(g)\pd{u^i_{l}(g)}{g^s}=Y^s_{jkl}u^i_s(g)
\end{equation}
Let us rewrite these differential equations as commutation relations.

In the tangent algebra $\Gamma$ of $G$ define the the ternary \emph{Yamaguti brackets} \cite{Yam63} $[\cdot,\cdot,\cdot]$ by
\begin{equation*}
[x,y,z]^i\=6Y^i_{jkl}x^jy^kz^l
\end{equation*}
Multiply (\ref{yam3}) by $6x^jy^kz^l$. Then we have
\begin{align*}
[x,y,z]
&=6(x,y,z)+2[[x,y],z]\\
&=[x,[y,z]]-[y,[x,z]]+[[x,y],z]
\end{align*}
Now from (\ref{pre-red}) it is easy to infer 
\begin{thm}[reductivity]
The infinitesimal translations of a local analytic Moufang loop satisfy the reductivity conditions
\begin{subequations}
\label{red}
\begin{align}
6[Y(x;y),L_z]&=L_{[x,y,z]}\\
6[Y(x;y),R_z]&=R_{[x,y,z]}\\
6[Y(x;y),M_z]&=M_{[x,y,z]}
\end{align}
\end{subequations}
\end{thm}

\begin{proof}
Commutation relation (\ref{red}a) is evident from  (\ref{pre-red}).  To get (\ref{red}b) one has to repeat the above calculations, i.e find the linear approximation of  (\ref{gle2yam}) with respect to $h$. Then  (\ref{red}c) easily follows by adding (\ref{red}a) and  (\ref{red}b).
\end{proof}

\section{Sagle-Yamaguti identity and hidden associativity}

Define the triality conjugated translations
\begin{equation*}
M^+\=L-R,\quad 
L^+\=R-M,\quad
R^+\=M-L
\end{equation*}
One can easily see the inverse conjugation:
\begin{equation*}
3M\=R^+-L^+,\quad 
3L\=M^+-R^+,\quad
3R\=L^+-M^+
\end{equation*}

\begin{thm}[reductivity]
The infinitesimal translations of a local analytic Moufang loop satisfy the reductivity conditions
\begin{subequations}
\label{red+}
\begin{align}
6[Y(x;y),L^+_z]&=L^+_{[x,y,z]}\\
6[Y(x;y),R^+_z]&=R^+_{[x,y,z]}\\
6[Y(x;y),M^+_z]&=M^+_{[x,y,z]}
\end{align}
\end{subequations}
\end{thm}

\begin{proof}
Evident corollary from formulae (\ref{red}).
\end{proof}

From \cite{Paal2} we know

\begin{prop}
Let $(S,T)$ be a Moufang-Mal'tsev pair. Then
\begin{subequations}
\label{lrm+}
\begin{align}
6Y(x;y)
&=[M^{+}_{x},M^{+}_{y}]+M^{+}_{[x,y]}\\
&=[R^{+}_{x},R^{+}_{y}]+R^{+}_{[x,y]}\\
&=[L^{+}_{x},L^{+}_{y}]+L^{+}_{[x,y]}
\end{align}
\end{subequations}
for all $x,y$ in $M$.
\end{prop}

\begin{thm}[hidden associativity]
The Yamagutian $Y$ of $G$ obey the commutation relations
\begin{equation}
6[Y(x;y),Y(z,w)]=Y([x,y,x],w)+Y(z;[x,y,w])
\end{equation}
if and only if the following Sagle-Yamaguti identity \cite{Yam62,Yam63} holds:
\begin{equation}
\label{sagle-yamaguti}
[x,y,[z,w]]=[[x,y,z],w]+[z,[x,y,w]]
\end{equation}
\end{thm}

\begin{proof}
We calculate the Lie bracket $[Y(x;y),Y(z,w)]$ from the Jacobi identity
\begin{equation}
\label{jacobi-temp}
[[Y(x;y),L^+_z],L^+_w]+[[L^+_z,L^+_w],Y(x;y]+[[L^+_w,Y(x;y),L^+_z]=0
\end{equation}
and formulae (\ref{lrm+}). We have
\begin{align*}
6[[Y(x;y),L^+_z],L^+_w]
&=[L^+_{[x,y,z]},L_w]\\
&=6Y([x,y,z];w)-L^+_{[[x,y,z],w]}\\
6[[L^+_z,L^+_w],Y(x;y]
&=36[Y(z;w),Y(x,y)]-6[L^+_{[z,w]},Y(x;y)]\\
&=36[Y(z;w),Y(x,y)]-L^+_{[x,y,[z,w]]}\\
6[[L^+_w,Y(x;y),L^+_z]
&=6Y(z;[x,y,w])-L^+_{[z,[x,y,w]]}
\end{align*}
By substituting these relations into (\ref{jacobi-temp}) we obtain
\begin{equation*}
36[Y(x;y),Y(z,w)]-6Y([x,y,x],w)-6Y(z;[x,y,w])
=L^+_{[x,y,[z,w]]-[[x,y,z],w]-[z,[x,y,w]]}
\end{equation*}
The latter relation has to be triality invariant. This mean that
\begin{subequations}
\label{l+r+m+}
\begin{align}
L^+_a
&=R^+_a=M^+_a\\
&=36[Y(x;y),Y(z,w)]-6Y([x,y,x],w)-6Y(z;[x,y,w])
\end{align}
\end{subequations}
where
\begin{equation*}
a=[x,y,[z,w]]-[[x,y,z],w]-[z,[x,y,w]]
\end{equation*}
But it easily follows from (\ref{l+r+m+}a) that
\begin{equation*}
L_a=R_a=M_a=0
\end{equation*}
from which in turn follows that $a=0$.
\end{proof}

\begin{rem}
A.~Sagle \cite{Sagle} and K.~Yamaguti proved \cite{Yam62} that the identity (\ref{sagle-yamaguti}) is equivalent to the Mal'tsev identity. In terms of Yamaguti \cite{Yam63} one can say that the Yamagutian $Y$ is a \emph{generalized representation} of the (tangent) Mal'tsev algebra $\Gamma$ of $G$.
\end{rem}

\section*{Acknowledgement} 

Research was in part supported by the Estonian Science Foundation, Grant 6912.

\bigskip\noindent
Department of Mathematics\\
Tallinn University of Technology\\
Ehitajate tee 5, 19086 Tallinn, Estonia\\ 
E-mail: eugen.paal@ttu.ee

\end{document}